\documentclass[12pt]{amsart}
\usepackage{graphicx}
\usepackage{amssymb}
\usepackage{amsfonts}
\usepackage{amsmath}
\usepackage{array}
\usepackage{rotating}

\def\FQSym{{\bf FQSym}}
\def\MQSym{{\bf MQSym}}
\def\PQSym{{\bf PQSym}}
\def\CQSym{{\bf CQSym}}
\def\SQSym{{\bf SQSym}}
\def\ev{{\rm Ev}}
\def\qrpark{{\bf q}}
\def\congru{\equiv}
\def\ssh{\Cup}
\def\saug{\uplus}
\def\sconc{\bullet}
\def\Std{{\rm Std}}
\def\Park{{\rm Park}}
\def\LPQ{{\bf LPQ}}
\def\NC{{\rm NC}}
\def\convol{{*}}

\def\park{{\bf a}}
\def\vp{{\rm vp}}

\def\<{\langle}
\def\>{\rangle}
\def\R{{\mathbb R}}
\def\C{\operatorname{\mathbb C}}

\def\N{\operatorname{\mathbb N}}

\def\F{{\bf F}}
\def\S{{\bf S}}
\def\T{{\bf T}}
\def\G{{\bf G}}
\def\M{{\bf M}}
\def\P{{\bf P}}
\def\V{{\bf V}}
\def\MA{{\mathcal M}}
\def\SG{{\mathfrak S}}

\def\Sym{{\bf Sym}}
\def\QSym{{\it QSym}}

\def\dim{{\rm dim}}

\def\ch{\operatorname {ch}}

\def\PF{{\rm PF}}
\def\PPF{{\rm PPF}}
\def\MM{{\mathcal M}}
\def\RR{{\bf R}}
\def\shuff#1#2{\mathbin{
\hbox{\vbox{ \hbox{\vrule \hskip#2 \vrule height#1 width 0pt
}%
\hrule}%
\vbox{ \hbox{\vrule \hskip#2 \vrule height#1 width 0pt
\vrule }%
\hrule}%
}}}

\def\shuf{{\mathchoice{\shuff{7pt}{3.5pt}}%
{\shuff{6pt}{3pt}}%
{\shuff{4pt}{2pt}}%
{\shuff{3pt}{1.5pt}}}}%
\def\shuffle{\,\shuf\,}

\def\Tabvrule{\vrule width-0.4pt}       
\def\Tabhrule{\hrule \hrule height-0.4pt} 
\def\Tabstrut{\vrule height2.2ex 
                     depth0.8ex  
                     width0ex    
\relax}

\def\PasCase#1{\omit%
            $\vcenter{\hbox {\vbox to 0.4pt{}}
               \hbox{\makebox[3ex]{\Tabstrut$#1$}}}%
               \Tabvrule$}
\def\PasCasePoint{\PasCase{\cdot}}
\def\DessinCarre#1{%
    \vcenter{\hbox{}\hrule
             \hbox{\vrule\makebox[3ex]{\Tabstrut$#1$}\vrule}\Tabhrule}%
             \Tabvrule}
\def\GenRuban#1{\vcenter{\halign{&$\DessinCarre{##}$\cr#1}}\egroup}

\def\sTabvrule{\vrule width-0.4pt}
\def\sTabhrule{\hrule \hrule height-0.4pt}
\def\sTabstrut{\vrule height1.6ex depth0.6ex width0ex \relax}
\def\sDessinCarre#1{%
    \vcenter{\hbox{}\hrule
             \hbox{\vrule\makebox[2.3ex]%
                  {\sTabstrut$\scriptstyle#1$}\vrule}\sTabhrule}%
             \sTabvrule}
\def\sGenRuban#1{\vcenter{\halign{&$\sDessinCarre{##}$\cr#1}}\egroup}

\def\ruban{%
  \bgroup
  \let\ =\omit
  \let\\=\cr
  \let\x=\times
  \let\.=\PasCasePoint
  \offinterlineskip
  \GenRuban}

\def\sruban{%
  \bgroup
  \let\ =\omit
  \let\x=\times
  \let\\=\cr
  \offinterlineskip
  \sGenRuban}

\title{A Hopf algebra of parking functions}

\author[J.-C.~Novelli and J.-Y.~Thibon]%
{Jean-Christophe Novelli and Jean-Yves Thibon}

\address[] {Institut Gaspard Monge, Universit\'e de Marne-la-Vall\'ee \\
5 Boulevard Descartes \\Champs-sur-Marne \\77454 Marne-la-Vall\'ee cedex 2 \\
FRANCE}
\email[Jean-Christophe Novelli]{novelli@univ-mlv.fr}
\email[Jean-Yves Thibon]{jyt@univ-mlv.fr} 
\date{}

\begin{document}

\begin{abstract}
If the moments of a probability measure on $\R$ are interpreted as
a specialization of complete homogeneous symmetric functions, its
free cumulants are, up to sign, the corresponding specializations
of a sequence of Schur positive symmetric functions $(f_n)$.
We prove that $(f_n)$ is the Frobenius characteristic of the natural
permutation representation of $\SG_n$ on the set of prime parking functions.
This observation leads us to the construction of a Hopf algebra of parking
functions, which we study in some detail.
\end{abstract}

\maketitle

\section{Introduction}

The free cumulants $R_n$ of a probability measure $\mu$ on $\R$
are defined (see {\it e.g.,}~\cite{Spei}) by means of the generating series
of its moments $M_n$
\begin{equation}
G_\mu(z) :=\int_\R\frac{\mu(dx)}{z-x}=z^{-1}+\sum_{n\ge 1}M_nz^{-n-1}
\end{equation}
as the coefficients of its compositional inverse 
\begin{equation}
K_\mu(z) :=G_\mu(z)^{\<-1\>}=z^{-1}+\sum_{n\ge 1}R_n z^{n-1}\,.
\end{equation}
It is in general instructive to interpret the coefficients of
a formal power series as the specializations of the elements of some 
generating family of the algebra of symmetric functions.
In this context, it is the interpretation
\begin{equation}
M_n=\phi(h_n)=h_n(A)
\end{equation}
which is relevant.
Indeed, the process of functional inversion (Lagrange
inversion) admits a simple expression within this formalism
(see~\cite{Mcd}, ex.~24 p.~35). If the symmetric functions
$h_n^*$ are defined by the equations
\begin{equation}
u=tH(t) \ \Longleftrightarrow \ t=uH^*(u)
\end{equation}
where $H(t):=\sum_{n\ge 0}h_nt^n$, $H^*(u):=\sum_{n\ge 0}h_n^*u^n$,
then, using the $\lambda$-ring notation,
\begin{equation}
h_n^*(X)=\frac{1}{n+1}(-1)^{n}e_n((n+1)X) :=\frac{1}{n+1} [t^n] E(-t)^{n+1}
\end{equation}
where $E(t)$ is defined by $E(t)H(t)=1$.
This defines an involution $f\mapsto f^*$ of the ring of symmetric
functions.

Now, if one sets $M_n=h_n(A)$ as above, then
\begin{equation}
G_\mu(z)=z^{-1}H(z^{-1})=u 
\ \Longleftrightarrow \ 
z= K_\mu(u)=\frac{1}{u}E^*(-u) =
u^{-1}+\sum_{n\ge 1}(-1)^ne_n^* u^{n-1}\,.
\end{equation}
Hence,
\begin{equation}  
R_n=(-1)^ne_n^*(A)\,.
\end{equation}
It follows immediately from the explicit formula (see~\cite{Mcd} p. 35)
\begin{equation}
-e_n^*=\frac1{n-1}\sum_{\lambda\vdash n}
\binom{n-1}{l(\lambda)} \binom{l(\lambda)}{m_1, m_2,\ldots,m_n} e_\lambda
\end{equation}
(where $\lambda=1^{m_1}2^{m_2}\cdots n^{m_n}$) that $-e_n^*$
is Schur positive. Clearly, $-e_n^*$ is the Frobenius characteristic
of a permutation representation $\Pi_n$, twisted by the sign character.
Let us set
\begin{equation}
(-1)^{(n-1)}R_n = -e_n^*=:\omega(f_n)
\end{equation}
so that $f_n$ is the character of $\Pi_n$.
We start with a construction of this representation in terms of parking
functions. This leads us to the definition of a Hopf algebra of parking
functions that generalizes the constructions of~\cite{MR,NCSF6}. We expect
that this combinatorics can be generalized to other root systems, at least for
type B (see, \emph{e.g.},~\cite{Bi1}).

\smallskip
We note that our construction of $\Pi_n$ is merely a variation about
previously known results (see in particular~\cite{Len,PP}). However, since
this is this precise version that led us to the Hopf algebra of parking
functions and some of its properties, we decided to present it in detail.

\smallskip
Although many definitions will be recalled, we shall assume that the reader is
familiar with the notation of~\cite{NCSF1,NCSF6}.

\medskip
{\footnotesize
{\it Acknowledgements.-}
This project has been partially supported by EC's IHRP Programme, grant
HPRN-CT-2001-00272, ``Algebraic Combinatorics in Europe".
The problem of constructing the representation $\Pi_n$ was suggested
by S. Kerov during his stay in Marne-la-Vall\'ee in 1996.
The question was forgotten for a long time without any attempt
of solution, and rediscovered recently on the occasion of talks by S.
Ferri\`eres and P. Biane.
Thanks also to P. Biane for providing the reference \cite{PP}.}

\section{Parking functions}


\subsection{Parking functions} 

A {\em parking function} on $[n]=\{1,2,\ldots,n\}$ is a word
$\park=a_1a_2\cdots a_n$ of length $n$ on $[n]$ whose nondecreasing
rearrangement $\park^\uparrow=a'_1a'_2\cdots a'_n$ satisfies $a'_i\le i$ for
all $i$.
Let $\PF_n$ be the set of such words.
It is well-known that $|\PF_n|=(n+1)^{n-1}$, and that the
permutation representation of $\SG_n$ naturally supported  by $\PF_n$
has Frobenius characteristic $(-1)^n\omega(h_n^*)$ (see~\cite{Haiman}).

\subsection{Prime parking functions} 

Gessel introduced in 1997 (see~\cite{Stan2}) the notion of
\emph{prime parking function}. 
One says that $\park$ has a {\em breakpoint} at $b$ if $|\{\park_i\le b\}|=b$.
Then, $\park\in \PF_n$ is said to be prime if its only breakpoint is $b=n$. 

Let $\PPF_n\subset\PF_n$ be the set of prime parking functions on $[n]$.
It can easily be shown that $|\PPF_n|=(n-1)^{n-1}$ (see~\cite{Stan2,Kal}).

\subsection{Operations on parking functions} 

For a word $w$ on the alphabet ${1,2,\ldots}$, denote by $w[k]$ the word
obtained by replacing each letter $i$ by $i+k$.
If $u$ and $v$ are two words, with $u$ of length $k$, one defines
the {\em shifted concatenation}
\begin{equation}
u\sconc v = u\cdot (v[k])
\end{equation}
and the {\em shifted shuffle}
\begin{equation}
u\ssh v= u\shuffle (v[k])\,.
\end{equation}
It is immediate to see that the set of permutations is closed under both
operations, and that the subalgebra spanned by those elements is isomorphic to
the convolution algebra of symmetric groups (see~\cite{MR}) or to Free
Quasi-Symmetric Functions (see~\cite{NCSF6}).

It is equally immediate to see that the set of all parking functions is
closed under these operations and that the prime parking functions exactly are
the parking functions that do not occur in any nontrivial shifted shuffle of
parking functions. These properties allow us to define a Hopf algebra of
parking functions (see Section~\ref{hopfalg}).

Let us now move to representation theory.

\subsection{The module of prime parking functions} 

Recall that the expression of complete symmetric functions in the
basis $e_\lambda$ is the commutative image of the formula
\begin{equation}
(-1)^nS_n =\sum_{I\vDash n}(-1)^{l(I)}\Lambda^I
\end{equation}
which, applied to $h_n^*$, gives
\begin{equation}
\label{caracpfn}
\ch(\PF_n)=(-1)^n\omega(h_n^*)=\sum_{I\vDash n}{f_{i_{1}}\cdot f_{i_{2}}\cdots
f_{i_r}}.
\end{equation}

Now, let us interpret this last formula.
Parking functions can be classified according to the factorization of their
nondecreasing reorderings $\park^\uparrow$ with respect to the operation of
shifted concatenation. That is, if
\begin{equation}
\park^\uparrow = w_1 \sconc w_2 \sconc\cdots\sconc w_r
\end{equation}
is the unique maximal factorization of $\park^\uparrow$, each $w_i$ is a
nondecreasing prime parking function.
Let us define $i_k=|w_k|$ and let $I=(i_1,\ldots,i_r)$. 
We shall say that $\park$ is of \emph{type} $I$ and denote by $\PPF_I$ the set
of parking functions of type $I$.

Then, the set $\PPF_{n}$ of prime parking functions of size $n$ obviously
is a sub-permutation representation of $\PF_n$, and it remains to compute its
Frobenius characteristic.
We prove that it is $f_n$, so that $\Pi_n$ can be identified with
$\PPF_n$.
It is sufficient to show that the number of prime parking
functions whose reordered evaluation is a given partition $\lambda$ is equal
to
$\frac{1}{n-1} \binom{n-1}{l(\lambda)} \binom{l(\lambda)}{m_1,m_2,\ldots,m_n}$
where $\lambda=1^{m_1}2^{m_2}\cdots n^{m_n}$.
Indeed, this number corresponds to the number of ways of putting the
$\lambda_i$ over $n-1$ places in a circle ; there is one circular word
associated with each circle whose reading is a prime parking function
(see~\cite{DM}). It then easily comes that
\begin{equation}
\ch(\PPF_n) = f_{n}\,,
\end{equation}
so that $\Pi_n$ can be identified with $\PPF_n$, as claimed before.

As a consequence, the set $\PPF_I$ of parking functions of type $I$ is a
sub-permutation representation of $\PF_n$ too, and its Frobenius
characteristic is
\begin{equation}
\ch(\PPF_I) = f_{i_1} \ldots f_{i_r}\,.
\end{equation}

Summing over all compositions $I$ of $n$ finally gives the right
interpretation of Equation~(\ref{caracpfn}).
A more transparent proof is given in Section~\ref{compos-sec}. 

\section{A Hopf algebra of parking functions} 
\label{hopfalg}

\subsection{The algebra $\PQSym$} 

We can embed the algebra of \emph{Free Quasi-Symmetric functions} $\FQSym$
of~\cite{NCSF6} inside the algebra spanned by the elements $\F_\park$
($\park\in\PF$), whose multiplication rule is defined by
\begin{equation}
\label{prodF}
\F_{\park'}\F_{\park''}:=\sum_{\park\in\park'\ssh\park''}\F_\park\,.
\end{equation}
We shall call this algebra $\PQSym$ (Parking Quasi-Symmetric functions).

For example,
\begin{equation}
\F_{12}\F_{11}= \F_{1233} + \F_{1323} + \F_{1332} + \F_{3123} + \F_{3132}
+ \F_{3312}\,.
\end{equation}
%

\subsection{The coalgebra $\PQSym$} 

There is a comultiplication on $\PQSym$ that naturally extends the
comultiplication of $\FQSym$. Recall (see~\cite{MR,NCSF6}) that if $\sigma$ is
a permutation,
\begin{equation}
\label{CoprodF}
\Delta\F_{\sigma} = \sum_{u\cdot v=\sigma}{\F_{\Std(u)} \otimes \F_{\Std(v)}},
\end{equation}
where $\Std$ denotes the usual notion of standardization of a word.

Given a word $w$, it is possible to define a notion of \emph{parkization}
$\Park(w)$, a parking function that coincides with $\Std(w)$ when $w$ is a
word without repetition.

For $w=w_1w_2\cdots w_n$ on $\{1,2,\ldots\}$, let us define
\begin{equation}
\label{dw}
d(w):=\min \{i | \#\{w_j\leq i\}<i \}\,.
\end{equation}
If $d(w)=n+1$, then $w$ is a parking function and the algorithm terminates,
returning~$w$. Otherwise, let $w'$ be the word obtained by decrementing all
the elements of $w$ greater than $d(w)$. Then $\Park(w):=\Park(w')$. Since
$w'$ is smaller than $w$ in the lexicographic order, the algorithm terminates
and always returns a parking function.

\smallskip
For example, let $w=(3,5,1,1,11,8,8,2)$. Then $d(w)=6$ and the word
$w'=(3,5,1,1,10,7,7,2)$.
Then $d(w')=6$ and $w''=(3,5,1,1,9,6,6,2)$. Finally, $d(w'')=8$ and
$w'''= (3,5,1,1,8,6,6,2)$, that is a parking function.
Thus, $\Park(w)=(3,5,1,1,8,6,6,2)$.

\smallskip
Now, the comultiplication on $\PQSym$ in defined as
\begin{equation}
\Delta \F_{\park}:= \sum_{u\cdot v=\park} \F_{\Park(u)} \otimes \F_{\Park(v)},
\end{equation}

For example,
\begin{equation}
\Delta\F_{3132} = 1\otimes\F_{3132} + \F_{1}\otimes\F_{132} +
\F_{21}\otimes\F_{21} + \F_{212}\otimes\F_{1} + \F_{3132}\otimes 1\,.
\end{equation}

One can easily check that the product and the comultiplication of $\PQSym$ are
compatible, so that $\PQSym$ is endowed with a bialgebra structure.

\subsection{The Hopf algebra $\PQSym$}

Since $\PQSym$ is endowed with a bialgebra structure naturally graded by the
size of parking functions, one defines the antipode as the inverse of the
identity for the convolution product and then endow $\PQSym$ with a Hopf
algebra structure.

The formula for the antipode can be written on the basis of $\F_{\park}$
functions, as
\begin{equation}
\nu (\F_{\park}) = \sum_{r ; u_1\cdots u_r=\park ; |u_i|\geq1} (-1)^r\,
\F_{\Park(u_1)}\F_{\Park(u_2)}\cdots \F_{\Park(u_r)}
\end{equation}

For example,
\begin{equation}
\nu (\F_{122}) = -\F_{122} + \F_{1}\F_{11} + \F_{12}\F_{1} - \F_{1}^3 =
\F_{212} + \F_{221} - \F_{213} - \F_{231} - \F_{321}\,.
\end{equation}

\subsection{The graded dual $\PQSym^*$}

Let $\G_{\park}=\F_{\park}^* \in\PQSym^*$ be the dual basis of $(\F_\park)$.
If $\langle\,,\,\rangle$ denotes the duality bracket, the product on
$\PQSym^*$ is given by
\begin{equation}
\label{prodG}
\G_{\park'} \G_{\park''} = \sum_{\park}
    \langle\, \G_{\park'}\otimes\G_{\park''}, \Delta\F_\park \,\rangle\,
    \G_\park
= \sum_{\park \in \park'\convol\park''} \G_\park\,,
\end{equation}
where the \emph{convolution} $\park'\convol\park''$ of two parking functions
is defined as
\begin{equation}
\park'\convol\park'' = \sum_{u,v ;
\park=u\cdot v, \Park(u)=\park', \Park(v)=\park''} \park\,.
\end{equation}

For example,
\begin{equation}
\begin{split}
\G_{12} \G_{11} &= \G_{1211} + \G_{1222} + \G_{1233} + \G_{1311} + \G_{1322}\\
&+ \G_{1411} + \G_{1422} + \G_{2311} + \G_{2411} + \G_{3411}\,.
\end{split}
\end{equation}

When restricted to permutations, it coincides with the convolution
of~\cite{Re,MR}.
Remark that in particular,
\begin{equation}
\G_1^n = \sum_{\park\in\PF_n} \G_\park\,.
\end{equation}

Using the duality bracket once more, one easily gets the formula for the
comultiplication of $\G_\park$ as
\begin{equation}
\Delta \G_\park := \sum_{u,v ; \park\in u\ssh v}
                   {\G_{\Park(u)} \otimes \G_{\Park(v)}}\,.
\end{equation}
There also exists a direct way to define the comultiplication of $\G_\park$
using the breakpoints of Gessel (see~\cite{Stan2}). In particular, the number
of terms in the coproduct is equal to the number of breakpoints of the parking
function plus one.

For example,
\begin{equation}
\begin{split}
\Delta \G_{41252} &= 1 \otimes\G_{41252} + \G_{1}\otimes\G_{3141} +
                   \G_{122}\otimes\G_{12} \\
                 &+ \G_{4122}\otimes\G_{1} + \G_{41252}\otimes1\,,
\end{split}
\end{equation}
whereas $41252$ has $4$ breakpoints : $1$, $3$, $4$, and $5$.

\subsection{Algebraic structure}
\label{alg-sec}

Let us say that a word $w$ over $\N^*$ is \emph{connected} if it cannot be
written as a shifted concatenation $w=u\sconc v$, and \emph{anti-connected} if
its mirror image $\overline{w}$ is connected.

Then, $\PQSym$ is free over the set
\begin{equation}
\left\{ \F_{\bf c}\, |\, {\bf c}\in\PF, \text{\ connected} \right\}
\end{equation}
ans $\PQSym^*$ is free over the set
\begin{equation}
\left\{ \G_{\bf d}\, |\, {\bf d}\in\PF, \text{\ anti-connected} \right\}
\end{equation}

This property proves that $\PQSym$ and $\PQSym^*$ are isomorphic as algebras. 
Moreover, it is possible to build an isomorphism $\varphi$ between $\PQSym$
and $\PQSym^*$ that is compatible with the product and the comultiplication.
So $\PQSym$ is isomorphic to $\PQSym^*$ as a \emph{Hopf algebra}.

When restricted to $\FQSym$, the isomorphism $\varphi$ is defined by
\begin{equation}
\label{phi-isom}
\varphi ( \F_\sigma) := \sum_{\park, \Std(\park)=\sigma^{-1}}{\G_\park}\,.
\end{equation}

\smallskip
The ordinary generating function for the numbers $c_n$ of connected parking
functions is
\begin{equation}
  \begin{split}
    \sum_{n\geq1} c_n t^n &= 1 -
         \left( \sum_{n\geq0}{(n+1)^{(n-1)} t^n} \right)^{-1} \\
    &= t + 2\, t^2 + 11\, t^3 + 92\, t^4 + 1014\, t^5 + 13795\, t^6  +
       223061\, t^7  + 4180785\, t^8 \\
    & + 89191196\, t^9  + 2135610879\, t^{10}   + 56749806356\, t^{11}   +
     1658094051392\, t^{12}\\
    & + O\left ({t}^{13}\right )\,.
  \end{split}
\end{equation}

\subsection{Multiplicative Bases}

Let $\park = {\bf a}_1 \sconc {\bf a}_2\sconc \cdots\sconc {\bf a}_r$ be the
maximal factorization of $\park$ into connected parking functions. We set
\begin{equation}
\F^{\park} = \F_{\park_1} \cdot \F_{\park_2} \cdots \F_{\park_r}\,,
\end{equation}
and
\begin{equation}
\G^{\overline{\park}} = \G_{\overline{\park_r}} \cdots
 \G_{\overline{\park_1}}\,.
\end{equation}
By a triangular argument, one can easily see that $(\F^{\park})$ (resp.
$(\G^{\overline{\park}})$), where $\park$ runs over the connected parking
functions, is a multiplicative basis of $\PQSym$ (resp. $\PQSym^*$).

\medskip
Now, if $\S_\park$ (resp. $\T_\park$) is the dual basis of $\F^{\park}$
(resp. $\G^{\overline{\park}}$) then
\begin{equation}
\{ \S_{\bf c} \,|\, {\bf c} \text{\ connected} \} \text{\ and\ }
\{ \T_{\bf c} \,|\, {\bf c} \text{\ connected} \}
\end{equation}
are bases of the primitive Lie algebras $\LPQ^*$ (resp. $\LPQ$)
of $\PQSym^*$ (resp. $\PQSym$). 

We conjecture, as in~\cite{NCSF6}, that both Lie algebras are free, on
generators whose degree generating function is
\begin{equation}
  \begin{split}
    1 - \prod_{n\geq1}{(1-t^n)}^{c_n} &= 1-(1-t)(1-t^2)^2(1-t^3)^{11} \cdots\\
    &= t + 2\,t^2 + 9\,t^3 + 80\,t^4 + 901\,t^5 + 12564\,t^6 + 206476\,t^7 \\
    & + 3918025\,t^8 + 84365187\,t^9 + 2034559143\,t^{10} +
      O\left ({t}^{11}\right )\,.
  \end{split}
\end{equation}

\subsection{Catalan Hopf algebra (non-crossing partitions)}

\subsubsection{The Hopf algebra $\CQSym$}

Parking functions are known to be related to non-crossing partitions
(see~\cite{Bi1,St,Stan2}).
There is a simple bijection between non-decreasing parking functions and
non-crossing partitions. Starting with a non-crossing partition, \emph{e.g.},
\begin{equation}
\pi = 13 | 2 | 45\,,
\end{equation}
one replaces all the letters of each block by its minimum, and reorders them
as a non-decreasing word
\begin{equation}
13 | 2 | 45 \to 11244
\end{equation}
which is a parking function. In the sequel, we identify non-decreasing parking
functions and non-crossing partitions via this bijection.

\smallskip
For a general $\park\in\PF_n$, let $\NC(\park)$ be the non-crossing partition
corresponding to $\park^\uparrow$ by the inverse bijection, \emph{e.g.},
$\NC(42141)=\pi$ as above.
Then, the elements of $\PQSym$
\begin{equation}
\P^\pi := \sum_{\park ; \NC(\park)=\pi} {\F_\park}
\end{equation}
span a sub-algebra of $\PQSym$,
isomorphic to the algebra of the free semigroup of non-crossing partitions
under the operation of concatenation of diagrams,
\begin{equation}
\P^{\pi'} \P^{\pi''} = \P^{\pi' \sconc \pi''}\,,
\end{equation}
that is equivalent to shifted concatenation on words.
Notice that $\P^\pi$ is the sum of all permutations of the
non-decreasing word corresponding to the given non-crossing partition.
We call this algebra the \emph{Catalan subalgebra} of $\PQSym$ and denote it
by $\CQSym$. The comultiplication is given on the basis $\P^\pi$ by
\begin{equation}
\label{coprodP}
\Delta\P^{\pi} = \sum_{u,v ; (u.v)^\uparrow=\pi}
{\P^{\Park(u)} \otimes \P^{\Park(v)}}\,,
\end{equation}
where $u$ and $v$ run over the set of non-decreasing words.

For example, one has
\begin{equation}
\begin{split}
\Delta\P^{1124} &= 1\otimes\P^{1124} +
        \P^1\otimes \left(\P^{112}+\P^{113}+\P^{123}\right) +
        \P^{11}\otimes\P^{12} \\
       & + \P^{12}\otimes\left(\P^{11}+2\P^{12}\right)
         + \left(\P^{112}+\P^{113}+\P^{123}\right)\otimes\P^1 +
        \P^{1124}\otimes1\,.
\end{split}
\end{equation}

One can easily check that the product and the comultiplication of $\CQSym$ are
compatible, so that $\CQSym$ is endowed with a graded bialgebra structure,
and therefore, with a Hopf algebra structure.
Formula~(\ref{coprodP}) immediately proves that the coalgebra $\CQSym$ is
co-commutative.


\subsubsection{The dual Hopf algebra $\CQSym^*$}

Let us denote by $\MM_\pi$ the dual basis of $\P^\pi$ in the
\emph{commutative} algebra $\CQSym^*$. Remark that $\CQSym^*$ is the quotient
of $\PQSym^*$ by the relations $\G_a\equiv\G_b$ if $a^\uparrow=b^\uparrow$. It
is then immediate (see Equation~(\ref{prodG})) that the multiplication is this
basis is given by
\begin{equation}
\label{prodMM}
\MM_{\pi'} \MM_{\pi''} = \sum_{\pi ;
\pi\in\pi'\convol\pi''}{\MM_{\pi^\uparrow}}\,.
\end{equation}

For example,
\begin{equation}
\begin{split}
\MM_{12} \MM_{11} &= \MM_{1112} + \MM_{1113} + \MM_{1114} + \MM_{1123} +
\MM_{1124}\\
&+ \MM_{1134} + \MM_{1222} + \MM_{1223} + \MM_{1224} + \MM_{1233}\,.
\end{split}
\end{equation}

This algebra can be embedded in the polynomial algebra $\C[x_1,x_2,\ldots]$ by
\begin{equation}
\MM_{\pi} = \sum_{\park(w)=\pi} {\underline{w}}\,,
\end{equation}
where $\underline{w}$ is the commutative image of $w$
(\emph{i.e.}, $i\mapsto x_i$).

For example,
\begin{equation}
\MM_{111} = \sum_{i} {x_i^3}\,.
\end{equation}
\begin{equation}
\MM_{112} = \sum_{i} {x_i^2x_{i+1}}\,.
\end{equation}
\begin{equation}
\MM_{113} = \sum_{i,j ; j\geq i+2} {x_i^2 x_j}\,.
\end{equation}
\begin{equation}
\MM_{122} = \sum_{i,j ; i<j} {x_i x_j^2}\,.
\end{equation}
\begin{equation}
\MM_{123} = \sum_{i,j,k; i<j<k} {x_ix_jx_k}\,.
\end{equation}

Notice that $\MM_{111}=M_3$; $\MM_{112}+\MM_{113}=M_{21}$ ; 
$\MM_{122}=M_{12}$ and $\MM_{123}=M_{111}$. In general, if
$\pi = \pi_1\sconc \cdots\sconc \pi_r$ is the factorization of $\pi$
in connected parking functions, let $i_k:=|\pi_k|$ and
$c(\pi):=(i_1,\cdots,i_k)$ a composition of $n$.
Then
\begin{equation}
\label{gammaEmb}
\gamma(M_I) := \sum_{c(\pi)=I} {\MM_{\pi}}
\end{equation}
gives an embedding of $\QSym$ into $\CQSym^*$.

\subsubsection{Catalan Ribbon functions}

In the classical case, the non-commutative complete fonctions split
into a sum of ribbon Schur functions, using a simple order on compositions.
To get an analogous construction in our case, we define a partial order on
non-decreasing parking functions.

Let $\pi$ be a non-decreasing parking function and $\ev(\pi)$ be its
evaluation vector. The successors of $\pi$ are the non-decreasing parking
functions whose evaluations are given by the following algorithm: given two
non-zero elements of $\ev(\pi)$ with only zeroes between them, replace the
left one by the sum of both and the right one by 0.

For example, the successors of $113346$ are $111146$, $113336$, and $113344$.

By transitive closure, the successor map gives rise to a partial order on
non-decreasing parking functions. We will write $\pi\preceq\pi'$ if $\pi'$ is
obtained from $\pi$ by successive applications of successor maps.

Now, define the Catalan Ribbon functions by
\begin{equation}
\label{catalRub}
\P^\pi =: \sum_{\pi'\succeq\pi} {\RR_{\pi'}}\,.
\end{equation}
This last equation completely defines the $\RR_{\pi}$.

\smallskip
The product of two $\RR$ functions is then
\begin{equation}
\RR_{\pi'} \RR_{\pi''} = \RR_{\pi'\sconc\pi''} +
\RR_{\pi'\triangleright\pi''}\,,
\end{equation}
where $\triangleright$ is the shifted concatenation defined by shifting all
elements of $\pi''$ by the difference between the greatest and the smallest
element of $\pi'$.

For example,
\begin{equation}
\RR_{11224} \RR_{113} = \RR_{11224668} + \RR_{11224446}\,.
\end{equation}

\subsection{Compositions}
\label{compos-sec}

Recall that non-crossing partitions can be classified according to the
factorization $\pi = \pi_1\sconc \cdots\sconc \pi_r$ into irreducible
non-crossing partitions.
We set
\begin{equation}
\V^I := \sum_{c(\pi)=I} \P^{\pi}
\end{equation}
as an element of $\PQSym$.
If one defines $\V_n = \V^{(n)}$, we have
\begin{equation}
\V_n = \sum_{\park\in\PPF_{n}} {\F_\park}
\end{equation}
and
\begin{equation}
\V^I = \V_{i_1} \cdots \V_{i_r} = \sum_{\park\in\PPF_I} {\F_\park}\,.
\end{equation}

At this point, it is useful to observe that if $C(w)$ denotes the descent
composition of a word $w$, the map
\begin{equation}
\eta: \F_\park \mapsto F_{C(\park)}\,,
\end{equation}
which is a Hopf algebra morphism $\PQSym\to\QSym$, maps $\V^I$ to the
Frobenius characteristic of the underlying permutation representation of
$\SG_n$ on $\PPF_I$.
\begin{equation}
\eta(\V^I) = \sum_{\park\in\PPF_I}\F_{C(\park)} = \text{ch}(\PPF_I)\,.
\end{equation}

As a consequence, the number of parking functions of type $I$ with descent set
$J$ is equal to the scalar product of symmetric functions
\begin{equation}
\< r_J, f^I\>
\end{equation}
where $f^I=f_{i_1}\cdots f_{i_r}=\text{ch}(\PPF_I)$ and $r_J$ is the ribbon
Schur function. This extends Prop.~3.2.(a) of~\cite{St}.
Remark that in particular,
\begin{equation}
\F_{\PF_{n}} := \sum_{\park\in\PF_{n}} {\F_\park} =
  \sum_{I\vDash n} {\V^I}\,,
\end{equation}
a realisation of Equation~(\ref{caracpfn}) as an identity in $\PQSym$.
By inversion, one obtains
\begin{equation}
\F_{\PPF_n} = \sum_{I\vDash n} {(-1)^{n-l(I)}\F_{\PF_I}}\,,
\end{equation}
where
\begin{equation}
\PF_I := \PF_{i_1} \ssh \PF_{i_2} \ssh \cdots \ssh \PF_{i_r}\,.
\end{equation}

These identities are easily visualized on the encoding of parking functions
with skew Young diagrams as in~\cite{PP} or in~\cite{HHLRU}.

The transpose $\gamma^*$ of the map $\gamma$ defined in
Equation~(\ref{gammaEmb}), is the map
\begin{equation}
\begin{split}
{\bf ch } :& \CQSym^*   \to \Sym \\
           & \P^\pi \mapsto  S^{c(\pi)}\,.
\end{split}
\end{equation}
which sends $\P^\pi$ to the characteristic non-commutative symmetric function
of the natural projective $H_n(0)$-module with basis
$\{\park\in\PF_n|NC(\park)=\pi\}$.

Then,
\begin{equation}
g := \sum_{n\geq0} g_n := \sum_{n\geq0} {\bf ch}(\F_{PF_n}) =
\sum_{I} {\bf ch}(V^I).
\end{equation}
is the series obtained by applying the non-commutative Lagrange inversion
formula of~\cite{Ges,PPR} to the generating series of complete functions,
\emph{i.e.}, $g$ is the unique solution of the equation
\begin{equation}
g = 1 + S_1g + S_2g^2 + \cdots = \sum_{n\geq0}{S_n g^n}\,.
\end{equation}

\subsection{Schr\"oder Hopf algebra (planar trees)}

Let $\congru$ denote the hypoplactic congruence (see~\cite{NCSF4,Nov}), and
denote by ${\sf P}(w)$ the hypoplactic $P$-symbol of a word $w$ (its
quasi-ribbon).  $P$-symbols of parking functions are called \emph{parking
quasi-ribbons}.

With a parking quasi-ribbon $\qrpark$, we associate the element
\begin{equation}
\P_{\qrpark} := \sum_{P(\park)=\qrpark} {\F_\park}\,.
\end{equation}
Then, the $\P_\qrpark$ form the basis of a Hopf sub-algebra of $\PQSym$,
denoted by $\SQSym$. Its dual $\SQSym^*$ is the quotient
$\PQSym/{\mathcal J}$ where $\mathcal J$ is the two-sided ideal generated by
\begin{equation}
\{\G_\park - \G_{\park'} | \park\congru\park' \}\,.
\end{equation}

If $\overline{\G_{\park}}$ denoted the equivalence class of $\G_\park$ modulo
$\mathcal J$, the dual basis of $(\P_\qrpark)$ is
\begin{equation}
{\bf Q}_{\qrpark} := \overline{\G_\park}\,,
\end{equation}
where $\park$ is any parking function such that $\park\equiv\qrpark$.

The dimension of the component of degree $n$ of $\SQSym$ and $\SQSym^*$ is the
little Schr\"oder number (or super-Catalan) $s_n$ : their Hilbert series is
\begin{equation}
\sum_{n\geq0} s_n t^n = \frac{1+t+\sqrt{1-6t+t^2}}{4t}
= 1 + t + 3t^2 + 11t^3 + 45t^4 + \cdots
\end{equation}

Indeed,
\begin{equation}
\begin{split}
\dim(\SQSym_n) &= \left\< \sum_{I\vDash n} F_I, {\bf ch}(\F_{\PF_n}) \right\>
= \left\< \frac{1}{2}\sum_{k=0}^n{e_kh_{n-k}}, \frac{1}{n+1}h_n((n+1)X)
\right\> \\
&= \frac{1}{2n+2} \sum_{k=0}^n{\binom{n+1}{k}\binom{2n-k}{n-k}} = s_n\,.
\end{split}
\end{equation}

The embedding of Formula~(\ref{phi-isom}) induces an embedding
\begin{equation}
\QSym \simeq \FQSym^*/({\mathcal J}\cap\FQSym^*) \rightarrow
\PQSym^*/{\mathcal J} = \SQSym^*\,.
\end{equation}
It is likely that $\SQSym$ is isomorphic to the free dendriform trialgebra
of~\cite{LRtri} as an algebra, but not as a coalgebra.

\subsection{$\PQSym^*$ as a combinatorial Hopf algebra}

Since $\FQSym$ can be embedded in $\PQSym$, we have a canonical Hopf embedding
of $\Sym$ in $\PQSym$ given by
\begin{equation}
S_n \mapsto \F_{12\cdots n}\,.
\end{equation}
With parking functions, we have other possibilities: for example,
\begin{equation}
j(S_n) := \F_{11\cdots 1}
\end{equation}
is a Hopf embedding, whose dual $j^*$ maps $\PQSym^*$ to $\QSym$ and
therefore endows $\PQSym^*$ with a different structure of combinatorial Hopf
algebra in the sense of~\cite{ABS}.

On the dual side, the transpose $\eta^*$ of the map $\eta$ defined in the
previous section corresponds to the Hopf embedding
\begin{equation}
S_n \mapsto \sum_{\Std(\park)=12\cdots n} \G_\park
\end{equation}
of $\Sym$ into $\PQSym^*$, which is therefore the restriction of the
self-duality isomorphism of formula~(\ref{phi-isom}) to the $\Sym$ subalgebra
$S_n=\F_{12\cdots n}$ of $\PQSym$.

\section{Realization of $\PQSym$}

It is possible to find a realization of $\PQSym$ in terms of $(0,1)$-matrices,
that is reminiscent of the construction of $\MQSym$ (see~\cite{Hiv,NCSF6}),
and that coincides with it when restricted to permutation matrices, providing
the natural embedding of $\FQSym$ in $\MQSym$.

Let $\MA_n$ be the vector space spanned by symbols $X_M$ where $M$ runs over
$(0,1)$-matrices with $n$ columns and an infinite number of rows, with $n$
nonzero entries, so that at most $n$ rows are nonzero.

Given such a matrix $M$, we define its \emph{vertical packing}
$P=\vp(M)$ as the finite matrix obtained by removing the null rows of $M$.

For a vertically packed matrix $P$, we define
\begin{equation}
\M_P = \sum_{\vp(M)=P} X_M\,.
\end{equation}

Now, given a $(0,1)$-matrix, we define its reading $r(M)$ as the word
obtained by reading its entries by rows, from left to right and top to bottom
and recording the numbers of the columns of the ones. For example, the reading
of the matrix
\begin{equation}
\begin{pmatrix}
0 & 1 & 1 & 0 \\
1 & 0 & 0 & 0 \\
0 & 1 & 0 & 0
\end{pmatrix}
\end{equation}
is $(2,3,1,2)$.

A matrix $M$ is said to be of \emph{parking type} if $r(M)$ is a parking
function.
Finally, for a parking function $\park$, we set
\begin{equation}
\F_{\park} := \sum_{r(P)=\park, \text{$P$ vertically packed}}{\M_P} =
\sum_{r(M)=\park} X_M\,.
\end{equation}

For example,
\begin{equation}
\F_{(1,2,2)} = \M_{
\begin{pmatrix}
1 & 1 & 0 \\
0 & 1 & 0 \\
\end{pmatrix}
}  + \M_{
\begin{pmatrix}
1 & 0 & 0 \\
0 & 1 & 0 \\
0 & 1 & 0 \\
\end{pmatrix}
}\,.
\end{equation}

The multiplication on $\MA=\oplus_{n} \MA_n$ is defined by columnwise
concatenation of the matrices:
\begin{equation}
X_M X_N = X_{M\cdot N}\,.
\end{equation}

In order to explicit the product of $\M_P$ by $\M_Q$, we first need a
definition.
Let $P$ and $Q$ be two vertically packed matrices with respective heights $p$
and $q$. The \emph{augmented shuffle} of $P$ and $Q$ is defined as follows:
let $r$ be an integer in $[\max(p,q),p+q]$. One inserts zero rows in $P$ and
$Q$ in all possible ways so that the resulting matrices have $p+q$ rows.
Let $R$ be the matrix obtained by concatenation of such pairs of matrices. The
augmented shuffle consists in the set of such matrices $R$ with nonzero rows.
We denote this set by $\saug(P,Q)$.

With this notation,
\begin{equation}
\M_P \M_Q = \sum_{R\in\saug(P,Q)} {\M_R}\,,
\end{equation}
and also
\begin{equation}
\F_{\park'}\F_{\park''}=\sum_{\park\in\park'\ssh\park''}\F_\park\,,
\end{equation}
that is the same as Equation~(\ref{prodF}).

\smallskip
Finally, concerning the comultiplication, one has first to define the
parkization $\Park(M)$ of a vertically packed matrix $M$, which consists in
iteratively removing column $d(r(M))$ until $M$ becomes a parking
matrix.

The comultiplication of a matrix $\M_P$ is then defined as:
\begin{equation}
\Delta \M_{P} = \sum_{Q\cdot R=P} \M_{\Park(Q)} \otimes \M_{\Park(R)}\,,
\end{equation}
It is then easy to check that
\begin{equation}
\Delta \F_{\park} =
  \sum_{u\cdot v=\park} \F_{\Park(u)} \otimes \F_{\Park(v)}\,,
\end{equation}
which is the same as Equation~(\ref{CoprodF}).

\subsection{Realization of $\FQSym$}

A parking matrix $M$ is said to be a \emph{word matrix} if there is exactly
one $1$ in each column.
Then $\FQSym$ is the Hopf subalgebra generated by the parking word matrices.


\footnotesize
\begin{thebibliography}{aa}
%
\bibitem{ABS}{\sc M. Aguiar, N. Bergeron}, and {\sc F. Sottile},
{\it Combinatorial Hopf Algebra and generalized Dehn-Sommerville relations\/},
preprint math.CO/0310016.
%
\bibitem{Bi1}{\sc P. Biane}, {\it Parking functions of types A and B},
Electronic J. Combin. {\bf 9} (2002), \# 7.
%
\bibitem{NCSF6}{\sc  G. Duchamp, F. Hivert}, and {\sc J.-Y. Thibon},
{\it Noncommutative symmetric functions VI: free quasi-symmetric functions and
related algebras},
Internat. J. Alg. Comput. {\bf 12} (2002), 671--717.
%
\bibitem{DM}{\sc A. Dvoretzky} and {\sc T. Motzkin},
{\it A problem of arrangements},
Duke Math. J. {\bf 14} (1947), 303--313.
%
\bibitem{NCSF1}{\sc I.M. Gelfand, D. Krob, A. Lascoux, B. Leclerc,
V.~S. Retakh}, and {\sc J.-Y. Thibon},
{\it Noncommutative symmetric functions},
Adv. in Math. {\bf 112} (1995), 218--348.
%
\bibitem{Ges}{\sc I. Gessel},
{\it Noncommutative Generalization and $q$-analog of the Lagrange Inversion
Formula},
Trans. Amer. Math. Soc. {\bf 257} (1980), no. 2, 455--482.
%
\bibitem{HHLRU}{\sc J. Haglund, M. Haiman, N. Loehr, J. B. Remmel},
and {\sc A. Ulyanov},
{\it A combinatorial formula for the character of the diagonal coinvariants},
preprint, 2003, math.CO/0310414.
%
\bibitem{Haiman}{\sc M. Haiman},
{\it Conjectures on the quotient ring by diagonal invariants},
J. Alg. Combin. {\bf 3} (1994), 17--76.
%
\bibitem{Hiv} {\sc F. Hivert},
{\it Combinatoire des fonctions quasi-sym\'etriques},
Th\`ese de Doctorat, Marne-La-Vall\'ee, 1999.
%
\bibitem{Kal}{\sc L. H. Kalikow}, {\it Enumeration of parking functions,
allowable permutation pairs, and labelled trees}, Ph. D. Thesis,
Brandeis University, 1999.
%
\bibitem{NCSF4}{\sc D. Krob} and {\sc J.-Y. Thibon},
{\it Noncommutative symmetric functions IV: Quantum linear groups and Hecke
algebras at $q=0$},
J. Alg. Comb. {\bf 6} (1997), 339--376.
%
\bibitem{Len}{\sc C.  Lenart},
{\it Lagrange inversion and Schur functions},
J. Algebraic Combin. {\bf 11} (2000), 1, 69--78.
%
\bibitem{LRtri}{\sc J.-L. Loday} and {\sc M.~O. Ronco},
{\it Trialgebras and families of polytopes},
preprint math.AT/0205043.
%
%
\bibitem{Mcd}{\sc I.G. Macdonald}, {\it Symmetric functions and Hall
polynomials}, 2nd ed., Oxford University Press, 1995.
%
\bibitem{MR}{\sc C. Malvenuto} and {\sc C. Reutenauer}, {\it
Duality between quasi-symmetric functions and the Solomon descent
algebra}, J. Algebra {\bf 177} (1995), 967--982.
%
\bibitem{Nov} {\sc J.-C. Novelli},
{\it On the hypoplactic monoid},
Disc. Math. {\bf 217}, 2000, 315--336.
%
\bibitem{PP}{\sc I. Pak} and {\sc A. Postnikov},
{\it Enumeration of trees and one amazing representation of $S_n$},
Proc. FPSAC'96 Conf., Minneapolis, MN, 1996, 385--389.
%
\bibitem{PPR}{\sc I. Pak, A. Postnikov}, and {\sc V.~S.  Retakh},
{\it Noncommutative Lagrange Theorem and Inversion Polynomials},
preprint, 1995, available at
{\tt http://www-math.mit.edu/${}^\sim$pak/research.html}.  
%
\bibitem{Re}{\sc C. Reutenauer},
{\it Free Lie algebras},
Oxford University Press, 1993.
%
\bibitem{Spei}{\sc R. Speicher},
{\it Multiplicative functions on the lattice on non-crossing partitions and
free convolution},
Math. Annalen {\bf 298} (1994), 141--159.
%
\bibitem{St}{\sc R. P. Stanley},
{\it Parking functions and noncrossing partitions},
Electronic J. Combin. {\bf 4} (1997), \# 2.
%
\bibitem{Stan2}{\sc R. P. Stanley}, {\it Enumerative combinatorics},
vol. 2, Cambridge University Press, 1999.
%
%
\end{thebibliography}
\end{document}